\documentclass[twocolumn]{svjour3}
\smartqed  
\pdfoutput=1
\usepackage{graphicx}
\usepackage{amssymb}
\usepackage{amsmath}
\usepackage{color}
\usepackage{hyperref}
\usepackage{mathptmx}
\journalname{}
\begin{document}

\title{Model-Free Continuation of Periodic Orbits in Certain Nonlinear Systems Using Continuous-Time Adaptive Control
\thanks{This work is supported by Agriculture and Food Research Initiative Competitive Grant no. 2014-67021-22109 from the USDA National Institute of Food and Agriculture.}
}

\titlerunning{Model-Free Continuation of Periodic Orbits in Certain Nonlinear Systems Using Continuous-Time Adaptive Control}        

\author{Yang Li        \and
        Harry Dankowicz 
}

\institute{Yang Li \and Harry Dankowicz \at Department of Mechanical Science and Engineering \\
              University of Illinois at Urbana-Champaign\\
              \email{danko@illinois.edu}
}

\date{Received: date / Accepted: date}

\maketitle

\begin{abstract}
This paper generalizes recent results by the authors on noninvasive model-reference adaptive control designs for control-based continuation of periodic orbits in periodically excited linear systems with matched uncertainties to a larger class of periodically excited nonlinear systems with matched uncertainties and known structure. A candidate adaptive feedback design is also proposed in the case of scalar problems with unmodeled nonlinearities. In the former case, rigorous analysis shows guaranteed performance bounds for the associated prediction and estimation errors. Together with an assumption of persistent excitation, there follows asymptotic convergence to periodic responses determined uniquely by an \emph{a priori} unknown periodic reference input and independent of initial conditions, as required by the control-based continuation paradigm. In particular, when the reference input equals the sought periodic response, the steady-state control input vanishes. Identical conclusions follow for the case of scalar dynamics with unmodeled nonlinearities, albeit with slow rates of convergence. Numerical simulations validate the theoretical predictions for individual parameter values. Integration with the software package \textsc{coco} demonstrate successful continuation along families of stable and unstable periodic orbits with a minimum of parameter tuning. The results expand the envelope of known noninvasive feedback strategies for use in experimental model validation and engineering design.

\keywords{Control-based continuation \and Model reference adaptive control \and Persistent excitation}
\end{abstract}

\section{Introduction}
\label{intro}

Control-based continuation provides a model-free approach for tracking periodic orbits of periodically excited nonlinear dynamical systems, independently of their orbital stability, under variations in experimentally accessible parameters \textcolor{black}{$p$} ~\cite{Barton201754,Barton2011,Barton2013,Bureau20145464,Kleyman2020,Lee2020261,Misra20082113,Renson20162775,Renson2019449,Schwartz1997664,Song2022,Tartaruga2019}. More advanced implementations also support tracking of special classes of periodic responses, e.g., those with vanishing phase lag relative to the excitation or corresponding to fold points in one-parameter bifurcation diagrams~\cite{Abeloos2022,Renson202183,Renson2017,Renson2016145,Sieber2008}. The approach embeds the experiment in a feedback control loop with control input $u$, parameterized by an experimentally accessible \emph{reference signal} $r$ and designed such that the response to zero control input ($u=0$) is that of the original system. Provided that the closed-loop dynamics (including in the control input) exhibit asymptotic convergence to limit cycle dynamics for initial conditions in some region $\mathcal{C}$, periodic reference signals in some set $\mathcal{R}$, and parameter values in some region $\mathcal{P}$, control-based continuation seeks to determine $r\in\mathcal{R}$ and $p\in\mathcal{P}$ such that $\lim_{t\rightarrow \infty}u(t)=0$ given an initial condition in $\mathcal{C}$. If such can be found, then the corresponding steady-state dynamics must coincide with a periodic orbit of the original system, albeit with orbital stability properties determined by the feedback control design. A feedback control design that supports such a determination is said to be \emph{non-invasive}, since it leaves the family of periodic orbits unchanged~\cite{Barton2012509,Bureau20135883,Sieber2009211}.

In an abstract sense, setting aside any concerns about accuracy and precision, the control-based continuation approach is thus straightforward to implement (but see~\cite{Beregi2021885,Renson20192811,Schilder2015251}). Firstly, formulate a non-invasive feedback control design, preferably with some \emph{a priori} understanding of how to interpret the relationship between the sought reference input and the desired, but unknown, limit cycle dynamics. Secondly, while maintaining dynamics in $\mathcal{C}$, perform iterative updates on $r$ and $p$ until the control input approaches $0$ asymptotically. Of course, in order to result in a finite-dimensional problem, the latter must be preceded by \emph{discretization} of the periodic reference input and a suitably chosen finite-time approximation of the periodic, steady-state control input.

In practice, the construction of a non-invasive feedback control design may necessitate some \emph{a priori} knowledge of the dynamics near the sought periodic orbits, particularly in order to ensure exponentially asymptotically stable limit cycle dynamics for the closed-loop system. This is the case for non-adaptive linear feedback control, for which the gains must be chosen to ensure that all associated Floquet multipliers lie inside the unit circle~\cite{Abeloos2022,Bureau20135883}. Additionally, particular feedback designs may fail to maintain dynamics in $\mathcal{C}$ or even guarantee bounded response of the closed-loop system with potentially disastrous consequences. The latter is true for linear feedback control. Finally, even as a non-invasive design may have been found for a particular parameter region $\mathcal{P}$, it may need to be retuned repeatedly to accommodate larger variations in $p$.

Once one moves beyond non-adaptive linear feedback control, a theoretical analysis is often restricted to particular classes of problems. In two recent papers~\cite{Li20202092,Li20212563}, the present authors investigated the use of adaptive feedback control to overcome the challenges outlined in the previous paragraph, specifically for tracking of fixed points in a class of single-input-single-output discrete-time dynamical systems and periodic orbits in a class of linear systems with matched uncertainties. As shown there, provable performance bounds were accompanied by a significant reduction in tuning effort. This came at the expense, however, of non-exponential rates of convergence, as well as a requirement that the frequency content of the reference input be sufficiently rich to result in persistent excitation of the closed-loop dynamics.

In this paper, we consider tracking of periodic orbits using adaptive feedback control for a larger class of dynamical systems than in our previous work, assuming nonlinearities of known structure in the case of problems of arbitrary dimension, and restricting attention to scalar problems in the case of uniformly bounded, unmodeled nonlinearities with uniformly bounded first-order partial derivatives. In both cases, we assume matched uncertainty, i.e., that an appropriately chosen control input could cancel the influence of parameter uncertainty. We rely on versions of the model-reference adaptive control approach and, in the case of nonlinearities of known structure, derive guaranteed performance bounds and demonstrate robustness to additive uniformly bounded disturbances. Integration with the \textsc{coco} continuation package shows successful tracking of stable and unstable periodic orbits with a minimum of manual tuning. In the case of unmodeled nonlinearities, it also highlights practical challenges associated with slow rates of convergence of the closed-loop dynamics.

The remainder of this paper is organized as follows. The class of nonlinear systems of initial interest is defined in Section~\ref{ch5:probFormu}, which also includes a discussion of a corresponding non-invasive, non-adaptive, linear feedback control design. Section~\ref{sec:adaptive control} describes a proposed non-invasive model-reference adaptive control algorithm and associated performance bounds. Numerical simulations in Section~\ref{sec:numeric} illustrate the performance of the controller at a fixed parameter value, while its use for control-based continuation is explored in Section~\ref{ch5:CBC} using an implementation in the \textsc{coco} software package~\cite{dankowicz2013recipes}. \textcolor{black}{Robustness of the control design under unmodeled, uniformly bounded, additive disturbances is considered in Section~\ref{sec:robust}.} For a class of scalar systems, Section~\ref{ch5:sysUnmodeled} relaxes the assumption that the structure of the nonlinearity be known to the control design and demonstrates the application of a proposed model-reference adaptive control design for parameter continuation. A brief concluding discussion follows in Section~\ref{ch5:concl}.

\section{Problem formulation}
\label{ch5:probFormu}

\subsection{\textcolor{black}{Model class}}
Following the discussion in \cite{Li20212563} for a class of linear systems, consider the dynamical system,
\begin{align}
\label{eq:SysDynNctrl}
\dot{q} = A q + b\left(\theta^\mathsf{T} Q(t,q) +\sigma\right),
\end{align}
where $A\in\mathbb{R}^{n\times n}$ is a known constant Hurwitz matrix, $b$ is a known constant vector, $\theta\in\mathbb{R}^m$ is an unknown constant vector, $\sigma$ is a known periodic function of period $T$, and $Q(t,q)$ represents a known \emph{nonlinear} function of $t$ and $q$ that is periodic in $t$ with period $T$. \textcolor{black}{The model form \eqref{eq:SysDynNctrl} reduces to that in \cite{Li20212563} when $Q(t,q)\equiv q$. For other choices of $Q(t,q)$,  \eqref{eq:SysDynNctrl} captures problems with a one-dimensional nonlinearity of arbitrary known complexity and with additional disturbance $\sigma$ along the same direction as the nonlinearity. In the case that $Q$ does not depend explicitly on $t$ (this is the case considered in the numerical example in Section~\ref{sec:numeric}), the response is driven by the known signal $\sigma$.}

\textcolor{black}{For $n=2$, 
\begin{equation}
    A=\begin{pmatrix}0 & 1\\-\omega_0^2 & -2\zeta\omega_0\end{pmatrix}\mbox{ and }b=\begin{pmatrix} 0 \\ 1\end{pmatrix}
\end{equation} \eqref{eq:SysDynNctrl} models a single-degree of freedom oscillator with natural frequency $\omega_0$ and damping constant $\zeta$ that is acted upon by an additional nonlinearity and periodic excitation with angular frequency $2\pi/T$. We may imagine the experimental determination of families of periodic responses under variations in $T$ as one goal of the control-based continuation analysis.}

\textcolor{black}{For mechanical systems with more than one degree of freedom, the form of \eqref{eq:SysDynNctrl} limits consideration to problems with only one source of nonlinearity and an excitation that is ``parallel'' to the nonlinearity. An example is given by the parametrically excited two-degree-of-freedom model \cite{Zaghari2019} obtained with 
\begin{equation}
    A=\begin{pmatrix}0 & 1 & 0 & 0\\
    -\frac{k_{\mathrm{01,lin}}+k_{12}}{m_1} & -\frac{c_{01}+c_{12}}{m_1} & \frac{k_{12}}{m_1} & \frac{c_{12}}{m_1}\\
    0 & 0 & 0 & 1\\
    \frac{k_{12}}{m_2} & \frac{c_{12}}{m_2} & -\frac{k_{02}+k_{12}}{m_2} & -\frac{c_{02}+c_{12}}{m_2}\end{pmatrix}
\end{equation}
in terms of the linear stiffness coefficients $k_{02}$, $k_{12}$, and $k_\mathrm{01,lin}$ and damping coefficients $c_{02}$, $c_{12}$, and $c_{01}$,
\begin{equation}
    b=\begin{pmatrix}0 &1 &0 &0\end{pmatrix}^\mathsf{T}\mbox{ and }\theta=-\frac{1}{m_1}\begin{pmatrix}k_\mathrm{01,nlin}\\k_\mathrm{PE,1,lin}\\k_\mathrm{PE,1,nlin}\end{pmatrix}
\end{equation}
in terms of the unknown stiffness coefficients $k_\mathrm{01,nlin}$, $k_\mathrm{PE,1,lin}$, and $k_\mathrm{PE,1,nlin}$,
\begin{equation}
Q(t,q)=\begin{pmatrix}q_1^3 & q_1\cos\Omega_\mathrm{PE} t & q_1^3\cos\Omega_\mathrm{PE}t\end{pmatrix}^\mathsf{T}    
\end{equation}
in terms of the excitation frequency $\Omega_\mathrm{PE}$, and $\sigma\equiv 0$. Here, $m_1$ and $m_2$ are two lumped masses along a clamped-free cantilever beam such that a harmonic current running through a coil imposes a time-varying, restoring force on the first mass that is nonlinear in displacement. Again, we may consider experimental continuation of periodic responses under variations in $\Omega_\mathrm{PE}$.}

\subsection{\textcolor{black}{Control objectives}}
Suppose that there exists a locally unique periodic (but \emph{a priori} unknown) solution $q^*$ of period $T$ \textcolor{black}{to \eqref{eq:SysDynNctrl}}. Due to the nonlinearity, $q^*(t)$ generally contains frequencies other than $\omega=2\pi/T$. The stability of $q^*$ is determined by the eigenvalues (the Floquet multipliers of $q^\ast$) of the monodromy matrix $\Phi(T)$, obtained from
\begin{align}
\dot{\Phi}=\left(A+b\theta^\mathsf{T} Q_q(t,q^\ast(t)))\right)\Phi,\quad \Phi(0)=\mathbb{I},
\end{align}
where the subscript $_q$ denotes the Jacobian with respect to $q$. As long as these eigenvalues lie inside the unit circle in the complex plane, then $q^\ast$ is locally asymptotically stable. \textcolor{black}{Asymptotic stability is not assumed, however, as we seek to use control-based continuation to locate and track such periodic solutions, even if unstable.}

Let $r(t)$ be a reference periodic function of period $T$ and define $x=q-r$. It follows that\textcolor{black}{
\begin{equation}
\label{eq:SysDynErrForm}
\dot{x} = A x + b\theta^\mathsf{T}\left(Q(t,x+r)-Q(t,r)\right)+ g,
\end{equation}
where
\begin{equation}
    g=-\dot{r}+Ar+b\left(\theta^\mathsf{T} Q(t,r)+\sigma\right)
\end{equation}
}is also periodic with period $T$ and identically equal to $0$ for $r\approx q^\ast$ if and only if $r=q^\ast$. By definition, the periodic function $x^\ast=q^\ast-r$ satisfies \eqref{eq:SysDynErrForm} and is locally asymptotically stable if all eigenvalues of $\Phi(T)$ are inside the unit circle in the complex plane. In the special case that $r=q^\ast$, it follows that $x^\ast(t)$ is identically equal to $0$.

We consider the introduction of a \textcolor{black}{matched} scalar control input $u$ as shown below,
\begin{align}
\label{eq:SysDynCtrl}
\dot{q} = A q + b\left(u + \theta^\mathsf{T} Q\left(t,q\right)+\sigma\right),
\end{align}
with the aim of having $u$ determined by $q$ and $r$, such that $u$ converges to a periodic steady-state signal that is uniquely determined by $r$ and equal to $0$ for $r\approx q^\ast$ and $q(0)\approx q^\ast(0)$ if and only if $r=q^\ast$, in which case $q(t)\rightarrow q^\ast(t)$ as $t\rightarrow\infty$. \textcolor{black}{We refer to such a control design as non-invasive along the sought periodic orbit.} By definition of $x$, we obtain
\begin{align}
\label{eq:NLErrSysCtrl}
\dot{x} = A x + b\theta^\mathsf{T}\left(Q(t,x+r)-Q(t,r)\right)+bu+ g
\end{align}
and it follows that the \textcolor{black}{construction of a non-invasive design along the sought periodic orbit} needs to ensure that $u(t)\rightarrow 0$ for $r\approx q^\ast$ and $x(0)\approx 0$ if and only if $g(t)\equiv 0$, and that $x(t)\rightarrow 0$ in this case.

\subsection{\textcolor{black}{Proportional feedback}}
As an example, \textcolor{black}{let}
\begin{equation}
\label{eq:proportional control}
u=-k^\mathsf{T} \left(Q(t,q)-Q(t,r)\right)
\end{equation}
for some to-be-determined constant vector $k$. Substitution yields
\begin{equation}
\dot{x}=Ax+ b(\theta-k)^\mathsf{T}\left(Q(t,x+r)-Q(t,r)\right)+g.
\end{equation}
It follows that $u(t)\rightarrow 0$ for $r\approx q^\ast$ and $x(0)\approx 0$ if $g(t)\equiv 0$ provided that all the eigenvalues of the monodromy matrix $\Phi(T)$ lie inside the unit circle, where $\Phi(t)$ in this case is governed by 
\begin{align}
\dot{\Phi}=\left(A+b(\theta-k)^\mathsf{T} Q_q(t,r)\right)\Phi,\quad \Phi(0)=\mathbb{I},
\end{align}
and that $x(t)\rightarrow 0$ in this case.  \textcolor{black}{The control design in \eqref{eq:proportional control}} is clearly \textcolor{black}{non-invasive along the periodic orbit in} the case when $k=\theta$. Importantly, given the local character of this control law, there is no \emph{a priori} degree of closeness that will guarantee the convergence of $x$ and $u$, nor are bounds available on transient deviations from $0$.

It is clear that $x(t)$ cannot converge to $0$ if $g\ne 0$ under the control law \eqref{eq:proportional control}. Under exceptional circumstances, it is still possible that $u(t)\rightarrow 0$ if $Q(t,x+r)-Q(t,r)$ converges to a signal in the orthogonal complement of $k$, in violation of our articulated objective. This possibility may be excluded on a case-by-case basis or, in the case of \eqref{eq:proportional control}, eliminated entirely by requiring that $x(t)\rightarrow 0$ if and only if $g(t)\equiv 0$.

\section{Model-reference adaptive control strategy}
\label{sec:adaptive control}

In the absence of knowledge about $\theta$, the selection of $k$ in \eqref{eq:proportional control} such that all Floquet multipliers have magnitude less than $1$ is, at best, trial-and-error. As an alternative, we consider a form of model-reference adaptive control~\cite{lavretsky2012robust} that relies on an adaptive estimate of $\theta$ to achieve the stated objective. \textcolor{black}{We show that this is non-invasive along the sought periodic orbit under generic conditions.} As a side benefit, we obtain guaranteed bounds on the transient and steady-state dynamics.

\subsection{\textcolor{black}{Control design}}
To this end, consider the control law
\begin{align}
\label{eq:nonlinearMRAC1ctrl}
u=-\hat{\theta}^\mathsf{T} \left(Q(t,q)-Q(t,r)\right),
\end{align}
where $\hat{\theta}(t)$ denotes a time-dependent estimate of $\theta$, such that
\begin{align}
\label{eq:nonlinearMRAC1adpt}
\dot{\hat{\theta}}&=-\Gamma e^\mathsf{T} PbQ(t,q),
\end{align}
defines the adaptive dynamics in terms of the \emph{adaptation gain} $\Gamma>0$. As usual, $P$ is a positive definite matrix that satisfies the algebraic Lyapunov equation $PA+A^\mathsf{T} P=-S$ for some positive definite matrix $S$. The \emph{prediction error} $e=x_m-x$ is defined in terms of the reference state $x_m$ governed by the differential equation
\begin{align}
\label{eq:nonlinearMRAC1Ref}
\dot{x}_m = Ax_m + b\tilde{\theta}^\mathsf{T} Q(t,r)+g,
\end{align}
where $\tilde{\theta}=\hat{\theta}-\theta$ is the \emph{estimation error}. It follows that
\begin{align}
\label{eq:nonlinearMRAC1matchErr}
\dot{e} = A e + b\tilde{\theta}^\mathsf{T} Q(t,q).
\end{align}
Notably, while $\tilde{\theta}$ appears explicitly in \eqref{eq:nonlinearMRAC1Ref}, terms involving $\theta$ cancel out of the sum of the last two terms, ensuring that \eqref{eq:nonlinearMRAC1adpt} and \eqref{eq:nonlinearMRAC1Ref} are implementable.

\subsection{\textcolor{black}{Lyapunov analysis}}
Now \textcolor{black}{let $\mathcal{B}$ denote a ball  of Euclidean radius $R$ such that $\theta\in\mathcal{B}$,  and assume} that the initial conditions are chosen so that $x_m(0)=x(0)$, i.e., $e(0)=0$, and $\hat{\theta}(0)\in\mathcal{B}$. Then, the Lyapunov function
\begin{align}
\label{eq:MRAClyap}
V=e^\mathsf{T} Pe+\frac{1}{\Gamma}\tilde{\theta}^\mathsf{T}\tilde{\theta}
\end{align}
satisfies 
\begin{align}
\label{eq:nonlinearMRAC1DV}
\dot{V}=-e^\mathsf{T} Se\leq 0,
\end{align}
and, consequently,
\begin{equation}
V(t)\le V(0)=\frac{1}{\Gamma}\|\tilde{\theta}(0)\|_2^2\le \frac{4R^2}{\Gamma}.
\end{equation}
It follows that
\begin{equation}
\|e(t)\|_2\le\frac{2R}{\sqrt{\lambda_{\mathrm{min}}(P)\Gamma}},
\end{equation}
where $\lambda_\mathrm{min}(P)$ is the smallest eigenvalue of $P$, and
\begin{equation}
\|\tilde{\theta}(t)\|_2\le 2R.
\end{equation}
Since $Q(t,r)$ and $g(t)$ are also uniformly bounded, it follows from \eqref{eq:nonlinearMRAC1Ref} and the fact that $A$ is Hurwitz that $x_m$ is bounded, and consequently, that $x$ and $q$ are bounded. Equation~\eqref{eq:nonlinearMRAC1ctrl} then implies that $u$ is bounded.

By the above analysis, $e$ and $\dot{e}$ are both bounded. This implies that $\ddot{V}$ is bounded and, by Barbalat's lemma~\cite{lavretsky2012robust}, that $\dot{V}(t)\rightarrow 0$, which in turn implies that $e(t)\rightarrow 0$, i.e., that $x_m(t)\rightarrow x(t)$ and, consequently, $\dot{\tilde{\theta}}(t)\rightarrow 0$. Moreover, since $\ddot{e}$ is bounded, it follows that $\dot{e}(t)\rightarrow 0$ and, consequently, $\tilde{\theta}^\mathsf{T}(t) Q(t,q)\rightarrow 0$.

\subsection{\textcolor{black}{Persistent excitation}}
\textcolor{black}{In order to conclude that $\tilde{\theta}(t)\rightarrow 0$, }suppose that the reference input $r$ is chosen so that the signal $Q(t,q)$ is \emph{persistently exciting}~\cite{Jenkins20182463,Narendra1987127}, i.e., that the smallest eigenvalue of the (at least positive semi-definite) matrix
\begin{equation}
\label{eq:persist}
\int_t^{t+T} Q(\tau,q(\tau))Q^\mathsf{T}(\tau,q(\tau))\,\mathrm{d}\tau
\end{equation}
is bounded from below by a positive number $\alpha$ for all $t$. \textcolor{black}{Although we cannot confirm the persistence of excitation of the signal $Q(t,q)$ \emph{a priori}, we may consider the integral  obtained by replacing $q(t)$ with $r(t)$ in \eqref{eq:persist}, since $r$ is assumed to be close to $q^\ast$ in practice. Since the integrand then becomes periodic, it suffices to compute its value for $t=0$.}

\textcolor{black}{It now follows from} the mean-value theorem and the observation that $\dot{\tilde{\theta}}(t)\rightarrow 0$ that \textcolor{black}{for every $\tau\in[t,t+T]$, there exists a $\sigma\in[t,t+T]$ such that $\tilde{\theta}(\tau)-\tilde{\theta}(t)=\tau\dot{\tilde{\theta}}(\sigma)\le TM(t)$ for some function $M(t)\rightarrow 0$ as $t\rightarrow\infty$. Consequently,}
\begin{align}
&\int_t^{t+T} \tilde{\theta}^\mathsf{T}(\tau) Q(\tau,q(\tau))Q^\mathsf{T}(\tau,q(\tau)) \tilde{\theta}(\tau)\,\mathrm{d}\tau\nonumber\\
&\quad=\tilde{\theta}^\mathsf{T}(t)\left(\int_t^{t+T} Q(\tau,q(\tau))Q^\mathsf{T}(\tau,q(\tau)) \,\mathrm{d}\tau\right)\tilde{\theta}(t)+\epsilon(t),
\end{align}
where $\epsilon(t)\rightarrow 0$ as $t\rightarrow \infty$. Consequently,
\begin{align}
\int_t^{t+T} &\tilde{\theta}^\mathsf{T}(\tau) Q(\tau,q(\tau))Q^\mathsf{T}(\tau,q(\tau)) \tilde{\theta}(\tau)\,\mathrm{d}\tau-\epsilon(t)\nonumber\\
&\qquad\geq\alpha\|\tilde{\theta}(t)\|_2^2.
\end{align}
Since the left-hand side converges to $0$ and the right-hand side is bounded from below by $0$, it follows that $\tilde{\theta}(t)\rightarrow 0$.

Since $\tilde{\theta}(t)\rightarrow 0$ and $\dot{\tilde{\theta}}$ is bounded, the classical result of Desoer \cite{desoer1969slowly} and Solo \cite{Solo1994331} applied to the governing equation
\begin{equation}
\dot{x}=A x-b\tilde{\theta}^\mathsf{T} \left(Q(t,x+r)-Q(t,r)\right)+g
\end{equation}
implies that $x(t)$ and, consequently, $u(t)$ both converge to periodic steady-state responses that are uniquely determined by $g$, and that $x(t)\rightarrow 0$ if and only if $g\equiv 0$. \textcolor{black}{We conclude that the model-reference adaptive control design is non-invasive along the sought periodic orbit.} Generically, we again expect that $u(t)\nrightarrow 0$ when $g\ne 0$.

\section{Numerical simulations}
\label{sec:numeric}

In this section, we explore the predictions from Section~\ref{sec:adaptive control} regarding the boundedness of the prediction and estimation errors, as well as the convergence of $\tilde{\theta}(t)$ to $0$ under suitable conditions on the reference input $r$.

As an example, consider the dynamical system \eqref{eq:SysDynNctrl} with
\begin{align}
A=\begin{pmatrix}
0 & 1 \\
-1.5 & -0.5
\end{pmatrix},~ b=\begin{pmatrix}
0 \\
1
\end{pmatrix},~
\theta=\begin{pmatrix}
0.5 \\
0.4 \\
-0.04
\end{pmatrix},
\end{align}
$Q(t,q)=\begin{pmatrix}q_1 & q_2 & q_1^3\end{pmatrix}^\mathsf{T}$, and $\sigma=\sin\omega t$. This corresponds to a harmonically excited Duffing oscillator in which the damping, stiffness, and nonlinearity coefficients are assumed unknown to the control design. \textcolor{black}{For $\omega=1$, there exists a periodic steady-state response $q^\ast$ given by} 
\begin{align}
q_1^*(t) &\approx -0.9928\cos{t}+2.9876\sin{t}+0.0336\cos{3t}\nonumber\\
&-0.0255\sin{3t}-0.0005\cos{5t}+0.00002\sin{5t},\label{eq:5th1}\\
q_2^*(t) &\approx 2.9876\cos{t}+0.9928\sin{t}-0.0765\cos{3t}\nonumber\\
&-0.1008\sin{3t}+0.0001\cos{5t}+0.0025\sin{5t}.\label{eq:5th2}
\end{align}
It is easy to check that $Q(t,q^\ast)$ is persistently exciting since the integral in \eqref{eq:persist} is independent of $t$ and positive definite with smallest eigenvalue \textcolor{black}{approximately} equal to $3.2$.

\begin{figure}[ht]
	\centering
	\includegraphics[width=.48\textwidth]{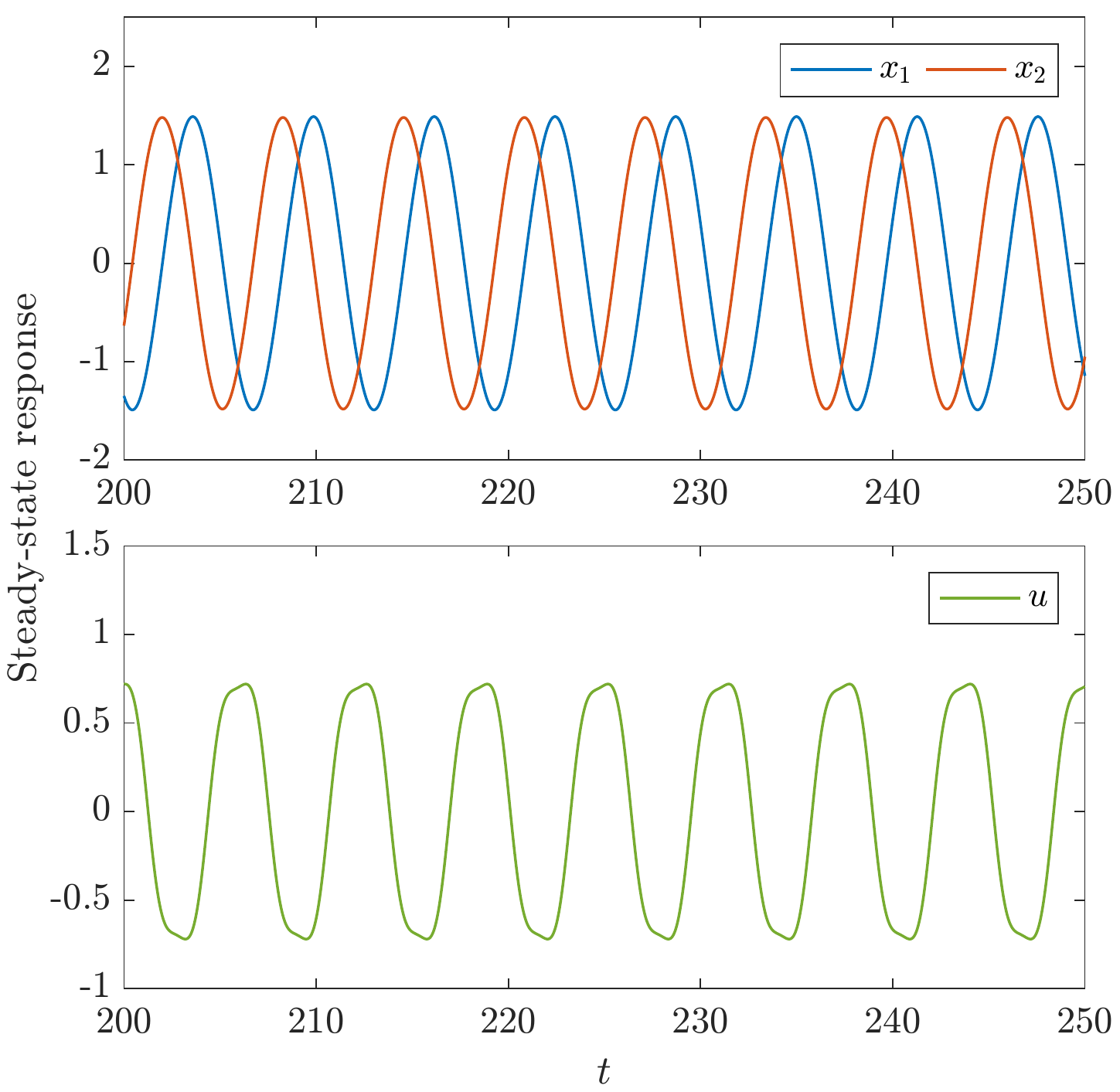}
	\caption{The steady-state responses of $x(t)$ and $u(t)$ under the proposed model-reference adaptive control strategy with $r(t)$ given in \eqref{eq:MRAC1DuffingNqstar}.}
	\label{fig:MRAC1DuffingssNqstar}
\end{figure}

Without loss of generality, we restrict attention to functions $r(t)$ chosen so that $\dot{r}-Ar$ is parallel to $b$ for all time, since this must be true of the desired reference input for which $g\equiv 0$. Consider the two choices of $r(t)=q^*(t)$ and 
\begin{equation}
\label{eq:MRAC1DuffingNqstar}
r(t)=\begin{pmatrix}\cos t+\sin t \\ \cos t-\sin t\end{pmatrix}.
\end{equation}
The steady-state solution to $\dot{x}=Ax+g$ is then given by $x(t)=0$ and the periodic function shown in Fig.~\ref{fig:MRAC1DuffingssNqstar}, respectively. The nonlinearity $Q(t,x+r)$ is persistently exciting also in the latter case, since the smallest positive eigenvalue of the integral \eqref{eq:persist} is again independent of $t$ and equal to $1.0$.

\begin{figure}[ht]
	\centering
	\includegraphics[width=.48\textwidth]{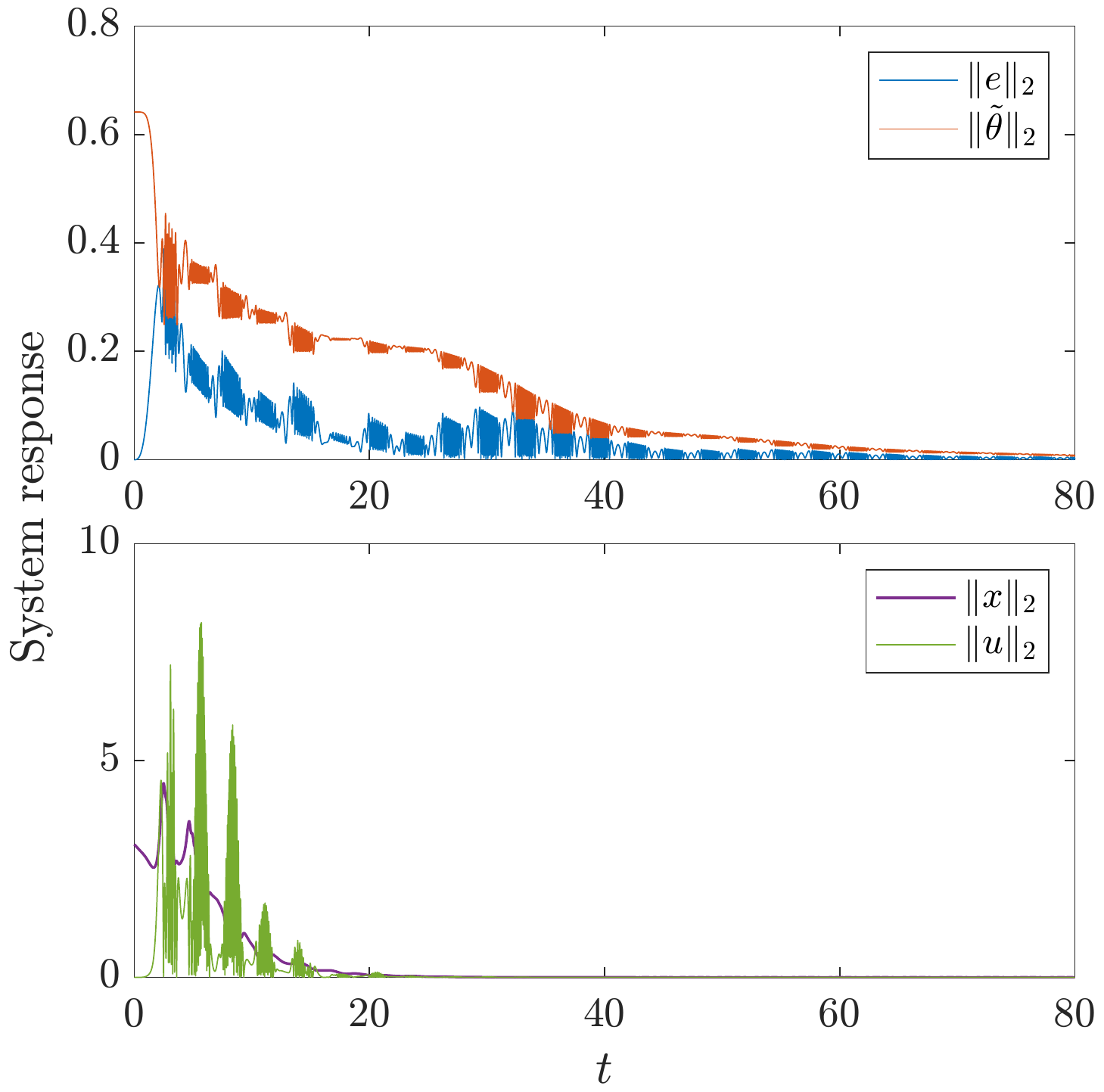}
	\caption{The system response under the proposed model-reference adaptive control strategy with $r(t)\approx q^*(t)$ given by the up-to-fifth-harmonic approximation in \eqref{eq:5th1}-\eqref{eq:5th2}. The predicted bounds on $\|e(t)\|_2$ and $\|\tilde{\theta}(t)\|_2$ approximately equal $1.28$ and $1.03$, respectively.}
	\label{fig:MRAC1Duffingqstar}
\end{figure}

Suppose that $q(0)=0$ and $\hat{\theta}(0)=0$, and let
\begin{align}
P=\begin{pmatrix}
8/3 & 1/3 \\
1/3 & 5/3
\end{pmatrix},
~ \Gamma=1.
\end{align}
The system response under the proposed model-reference adaptive control strategy with $r(t)= q^\ast(t)$ up to the fifth harmonic is shown in Fig.~\ref{fig:MRAC1Duffingqstar}. It is seen that $\|e(t)\|$ and $\|\tilde{\theta}(t)\|$ both go to $0$ as $t\rightarrow\infty$. Similarly, to within the resolution of the first five harmonics, $\|x(t)\|$ and $\|u(t)\|$ also converge to $0$. 

\begin{figure}[ht]
	\centering
	\includegraphics[width=.48\textwidth]{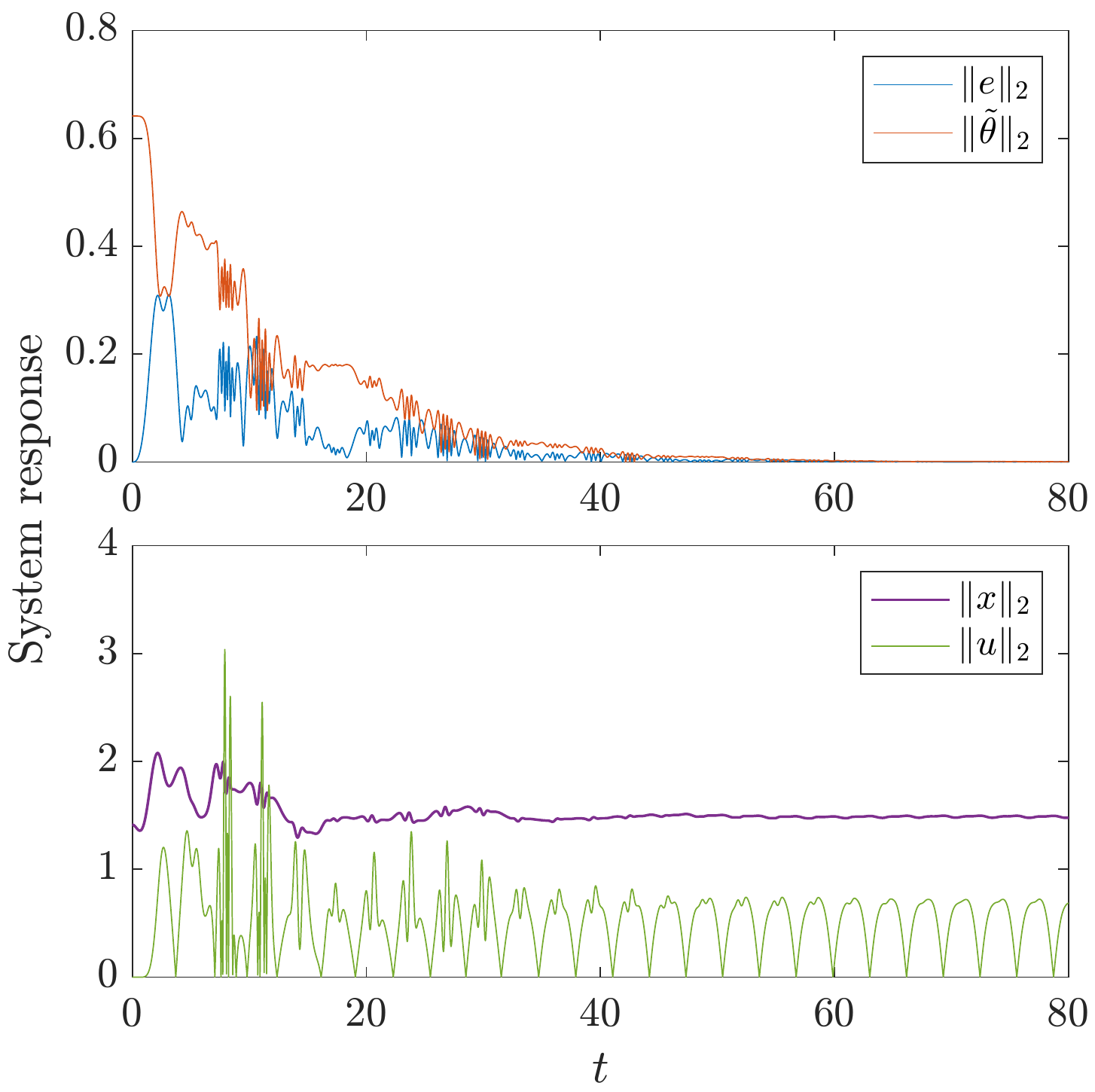}
	\caption{The system response under the proposed model-reference adaptive control strategy with $r(t)$ given in \eqref{eq:MRAC1DuffingNqstar}. The predicted bounds on $\|e(t)\|_2$ and $\|\tilde{\theta}(t)\|_2$ approximately equal $1.28$ and $1.03$, respectively.}
	\label{fig:MRAC1DuffingNqstar}
\end{figure}

The system response with $r(t)$ given in \eqref{eq:MRAC1DuffingNqstar} is shown in Fig.~\ref{fig:MRAC1DuffingNqstar}. It is seen that again, $\|e(t)\|$ and $\|\tilde{\theta}(t)\|$ both converge to $0$ as $t\rightarrow\infty$. Since $r$ deviates from a periodic steady-state response of the system, $\|x(t)\|$ and $\|u(t)\|$ converge to nonzero periodic responses, as predicted by the analysis in the previous section.

\begin{figure}[ht]
	\centering
	\includegraphics[width=.48\textwidth]{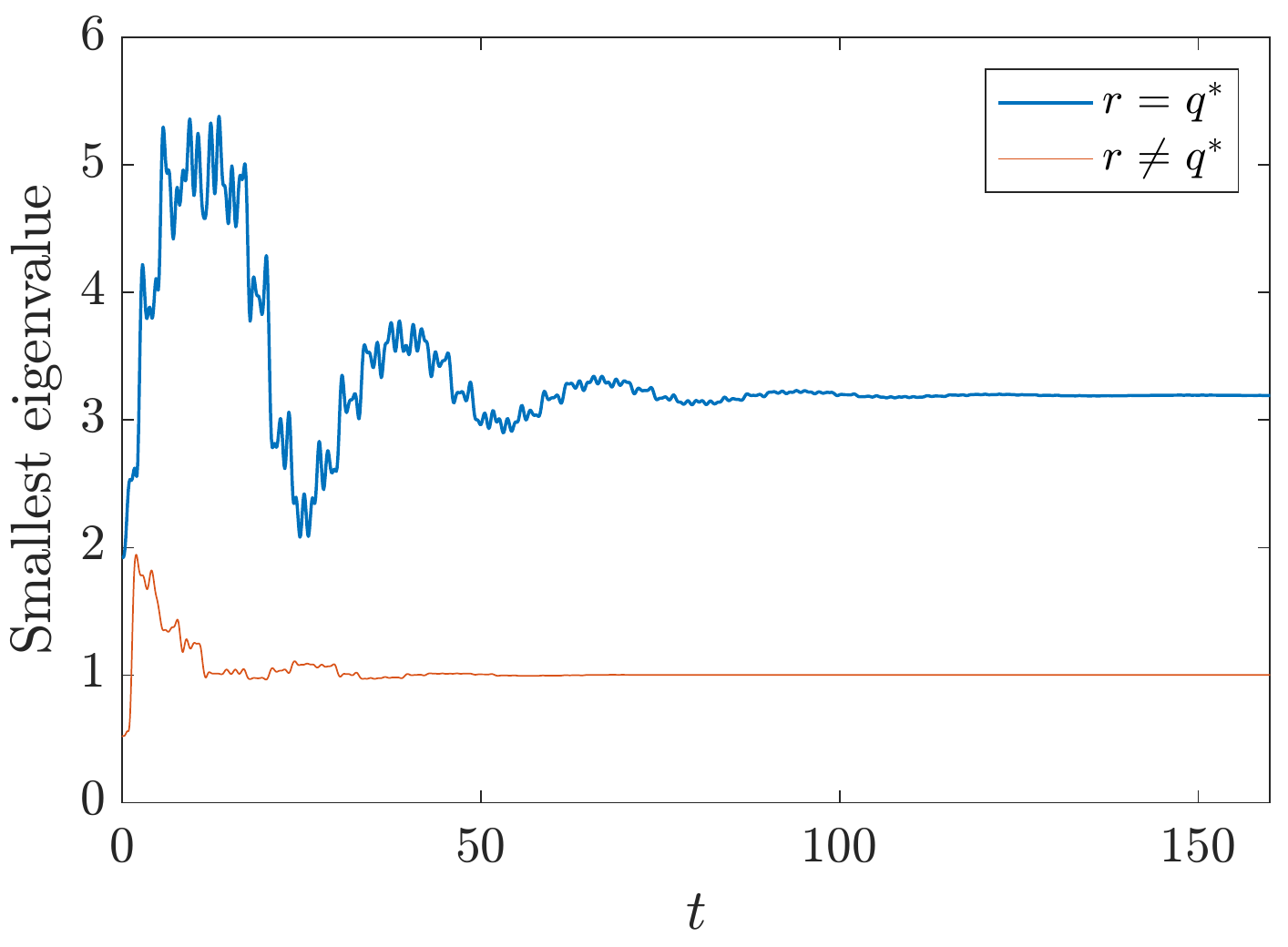}
	\caption{The smallest eigenvalue of the integral in \eqref{eq:persist} as a function of $t$, with $r(t)\approx q^*(t)$ given by the up-to-fifth-harmonic approximation in \eqref{eq:5th1}-\eqref{eq:5th2} and \eqref{eq:MRAC1DuffingNqstar}, respectively, under the model reference adaptive control.}
	\label{fig:alphaPE}
\end{figure}

For each of these two cases, Fig.~\ref{fig:alphaPE} shows the time histories of the smallest eigenvalue of the integral in \eqref{eq:persist}. Since these are bounded below by some positive numbers (and, indeed, converge to the predicted values obtained by substituting $q=x+r$ with $x$ being the unique solution to $\dot{x}=Ax+g$), $Q(t,q)$ is persistently exciting in both cases. This is also consistent with the general theory that ensures persistent excitation provided that the frequency content exceeds a multiple of the system dimension. Here, this is guaranteed by the nonlinearity. 

\section{Control-based continuation}
\label{ch5:CBC}

The theoretical treatment in Section~\ref{sec:adaptive control} shows that we may identify an \emph{a priori} unknown periodic response $q^\ast$ by the fact that $u(t)\rightarrow 0$ as $t\rightarrow\infty$ provided that the reference input $r$ happens to equal $q^\ast$. In this section, we use Newton's method to iteratively improve upon the reference input in order for the steady-state control input to fall within a threshold distance from $0$, thereby obtaining an approximation of $q^\ast$. As in the previous section, we restrict attention to $r$ chosen so that $\dot{r}-Ar$ is parallel to $b$, i.e., such that $r$ may be parameterized by a scalar periodic function of the same dimension as the control input. We proceed to use simulations of the closed-loop dynamics to estimate the coefficients of a truncated Fourier series of the steady-state control input and their derivatives with respect to the corresponding coefficients of $r$ and modify the coefficients of $r$ accordingly.

We combine the application of Newton's method with a pseudo-arclength continuation algorithm~\cite{dankowicz2013recipes} in order to trace $q^\ast$ under variations in a model parameter, say $\omega$ in the example in Section~\ref{sec:numeric}\textcolor{black}{, also past geometric folds where the assumption of existence and local uniqueness of $q^\ast$ for fixed $\omega$ would fail}. Here, an approximation to the tangent line of the graph of $(\omega,q^\ast)$ at a particular point on this graph is used to construct a \emph{predictor} $(\omega,r)$ some distance $h$ from $(\omega,q^\ast)$ along the tangent line. We proceed to require that all subsequent iterates of Newton's method lie on a line perpendicular to the tangent and intersecting the tangent at the predictor. We initialize the overall algorithm at some point $(\omega,q^\ast)$ obtained for example using forward simulation in the case that the corresponding $q^\ast$ is asymptotically stable. Importantly, convergence of the Newton iterations is independent of the open-loop stability of $q^\ast$. Also, since $r$ remains close to $q^\ast$, $g$ is close to $0$ for all iterations and the adaptive control strategy ensures the desired convergence.

Consistent with the implementation of the continuation algorithm in a physical experiment, we assume no direct control over the state $q$, for example its value at any moment in time. Other than the first simulation of the closed-loop dynamics, we initialize $q$ at its terminal value in the preceding simulation. We similarly initialize $\hat{\theta}$ for each simulation at its terminal value in the preceding simulation, thus ensuring close tracking by $x$ of the reference state $x_m$, since $\theta$ is assumed to vary smoothly with the model parameter.

Since it is not possible to obtain an exact match of $r$ and $q^\ast$, given the presence of harmonics of all orders, we select a truncation order that (empirically) yields sufficient information about $q^\ast$ while also ensuring that $Q(t,q)$ is persistently exciting. We avoid aliasing by applying the discrete Fourier transform to a fine sample of a period of the steady-state control input.

\begin{figure}[h]
  \includegraphics[width=0.48\textwidth]{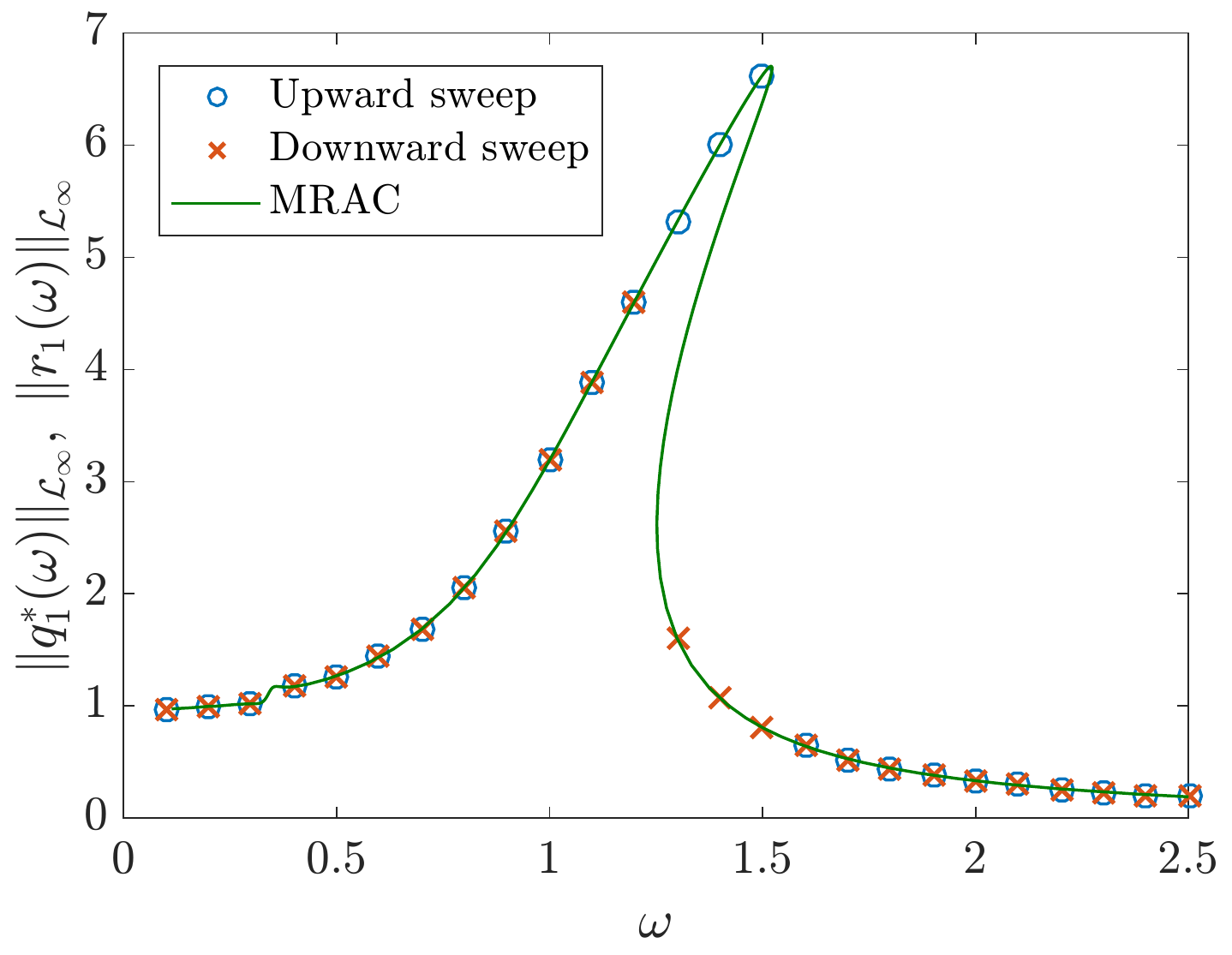}
\caption{A comparison between branches of periodic solutions of the harmonically excited nonlinear Duffing oscillator obtained using upward and downward sweeps in the excitation frequency $\omega$ (represented by $\|q^\ast_1(\omega)\|_{\mathcal{L}_\infty}$), and those obtained using the control-based continuation algorithm with the proposed model reference adaptive control strategy (represented by $\|r_1(\omega)\|_{\mathcal{L}_\infty}$).}
\label{fig:DuffingCBC}
\end{figure}

Figure~\ref{fig:DuffingCBC} shows the successful application of the control-based continuation algorithm, implemented in the software package \textsc{coco}~\cite{dankowicz2013recipes}, to the example in Section~\ref{sec:numeric} under variations in $\omega$ and with control parameters $P$ and $\Gamma$ as given there. Throughout continuation, we approximate $u$ by its truncated Fourier series up to the fifth harmonic and iteratively update the corresponding coefficients of $r$. In each iteration, and when approximating derivatives with respect to the Fourier coefficients of $r$, the Matlab integrator \texttt{ode45} (with relative tolerance $10^{-8}$ and absolute tolerance $10^{-10}$) is used to simulate the closed-loop transient dynamics for 10 periods, followed by a sampling of the control input during one additional period of simulation. We consider $r$ to have converged to $q^\ast$ when the norm of the Fourier coefficients of the steady-state control input is smaller than $10^{-6}$. The step size along the graph $(\omega,q^\ast)$ is adaptively determined by \textsc{coco} using default settings. As seen in the figure, the algorithm is able to trace out the solution branch independent of the open-loop stability of the periodic solutions (which are unstable along the middle branch in the range of coexisting periodic solutions).

\textcolor{black}{
\section{Robustness}
\label{sec:robust}
In practice, the presence of unmodeled disturbances may limit the utility of the proposed methodology, both in terms of the expected convergence of $u(t)$ to $0$ as $t\rightarrow \infty$ when the reference signal $r$ is chosen appropriately, and in terms of any guarantees on a bounded response.}

\textcolor{black}{Consider, for example, the introduction in the closed-loop dynamics of an additive, uniformly bounded, unknown disturbance $h(t,q)$:
\begin{align}
\label{eq:SysDynCtrlNoise}
\dot{q} = A q + b\left(u+\theta^\mathsf{T} Q\left(t, q\right)+\sigma\right)+h(t,q),
\end{align}
where $u$ is given in \eqref{eq:nonlinearMRAC1ctrl}. With $x=q-r$, $e=x_m-x$, and $x_m$ again governed by
\begin{equation}
    \dot{x}_m=Ax_m+b\tilde{\theta}^\mathsf{T}Q(t,r)-\dot{r}+Ar+b\left(\theta^\mathsf{T}Q(t,r)+\sigma\right),
\end{equation}
it then follows that
\begin{align}
\dot{e} = A e + b\tilde{\theta}^\mathsf{T} Q\left(t, q\right)-h(t,q).
\end{align}
Given the Lyapunov function in \eqref{eq:MRAClyap}, the adaptation law \eqref{eq:nonlinearMRAC1adpt} implies that
\begin{align}
\label{eq:appVdot}
\dot{V}=-e^\mathsf{T} Se+2e^\mathsf{T} Ph(t,q).
\end{align}
Because of the second term on the right-hand side, we can no longer claim that $\dot{V}\leq 0$ or $e\rightarrow 0$ as $t\rightarrow \infty$. It is still the case that $e$ is bounded, however. Indeed, if $h_b$ denotes an upper bound for $\|h(t,q)\|$, then
\begin{align}
\dot{V}(t)\leq -\lambda_{\mathrm{min}}(S)\|e\|^2+2\|e\|\lambda_{\mathrm{max}}(P)h_b.
\end{align}
Consequently, $\dot{V}<0$ when 
\begin{align}
\label{app:ineq}
\|e\|> 2\frac{\lambda_{\mathrm{max}}(P)}{\lambda_{\mathrm{min}}(S)}h_b.
\end{align}
While this implies that $\|e\|$ is eventually upper bounded by the right hand side of \eqref{app:ineq}, boundedness of $\hat{\theta}$ does not follow. To achieve this, and by implication the boundedness of $x$, $q$, and $u$, we may modify the adaptation law using a projection operator~\cite{lavretsky2012robust}
\begin{align}
\dot{\hat{\theta}}=\Gamma\,\mathrm{Proj}_\mathcal{B}\left(\hat{\theta}, -e^\mathsf{T} PbQ\left(t, q\right)\right),\, \hat{\theta}(0)=\hat{\theta}_0.
\end{align}
With this modification, the proposed model-reference adaptive control design guarantees a bounded response given a uniformly bounded additive disturbance $h(t,q)$, a reassuring prediction for any actual implementation in a physical experiment.}

\textcolor{black}{In the special case that $h(t,q)$ is periodic in $t$ with period $T$, the unique, locally attractive limit cycle dynamics of the closed-loop system obtained when $h\equiv 0$ and $r\approx q^\ast$ persist for sufficiently small $\|h\|_\infty$. In this case, if $r$ is chosen so that
\begin{equation}
    \dot{r}=Ar+b\left(\theta^\mathsf{T}Q(t,r)+\sigma\right)+h(t,r)
\end{equation}
then
\begin{align}
    \dot{x}=Ax&-b\tilde{\theta}^\mathsf{T}\left(Q(t,x+r)-Q(t,r)\right)\nonumber\\
    &+h(t,x+r)-h(t,r)
\end{align}
and we conclude that $x(t)\equiv 0$ along the perturbed limit cycle and $q\rightarrow r$ locally if and only if $r$ is chosen in this way.
}

\textcolor{black}{In the case that $h$ is not periodic in $t$, we cannot expect local persistence. In this case, the control-based continuation algorithm fails to trace out the solution branch, since it is no longer the case that $u(t)$ converges to $0$ as $t\rightarrow\infty$ when $r=q^*$ if the disturbance persists. Given additional information about the disturbance, it may be possible to substitute the requirement that $u(t)\rightarrow0$ with a condition on a suitably filtered version of $u(t)$. In either case, bounded performance is guaranteed.
}

\section{Systems with unmodeled nonlinearities}
\label{ch5:sysUnmodeled}

In this section, we attempt to relax the expectation that the form of the nonlinearity be known to the control design \textcolor{black}{and accessible to the feedback law}.  As a first step in this direction, consider the example system \textcolor{black}{
\begin{align}
\label{eq:scalarSys1}
\dot{q} = -q + \sin q + \sin\omega t+u,\,q\in\mathbb{R}.
\end{align}
When $u=0$, there exists a} periodic solution $q^\ast(t)$ for $\omega=1$ \textcolor{black}{given by}
\begin{align}
\label{eq:qstarapprox}
q^*(t) &\approx -0.9849\cos{t}+0.1160\sin{t}+0.0053\cos{3t}\nonumber\\
&+0.0115\sin{3t}+0.0002\cos{5t}-0.0003\sin{5t}.
\end{align}


Next, let $x=q-r$ and $u=-\hat{k}x$, and consider the closed-loop dynamics of the system 
\begin{align}
\label{eq:scalarSys2}
\dot{x} = -(1+\hat{k})x +g,\, \dot{\hat{k}}=\Gamma x^2
\end{align}
where \textcolor{black}{
\begin{equation}
g=\sin q+\sin \omega t-\dot{r}-r    
\end{equation}
}for some periodic reference input $r$. We note that  \textcolor{black}{$x\equiv0$ is a solution of this system provided that $g|_{q=r}\equiv0$}, i.e., that $r$ is a periodic solution of the open-loop dynamics. Indeed, when this is not the case, $\hat{k}$ is non-decreasing and must grow beyond all bounds as $t\rightarrow\infty$. Given the Lyapunov function
\begin{align}
V=x^2+\frac{1}{\Gamma}\hat{k}^2,
\end{align}
it follows that
\begin{align}
\dot{V}(t)&=-2x^2+2xg.
\end{align}
Let $g_b$ denote an upper bound on the magnitude of $g$. Then, since $\hat{k}$ is a non-decreasing function of time, \textcolor{black}{
\begin{align}
\label{eq:MRAC2xbd}
|x(t)|\leq \max\{g_b, |x(0)|\}+|\hat{k}(0)|/\sqrt{\Gamma}.
\end{align}
Indeed, for $\hat{k}(0)>0$, $|x(t)|$ cannot exceed $\max\{g_b, |x(0)|\}$, since $|x(t^*)|>g_b$ at some instant $t^*$ implies that $\dot{V}(t^*)<0$ and, consequently, that $x(t^*)\dot{x}(t^*)<0$. If, instead, $\hat{k}(0)<0$, then $V$ cannot exceed $\max\{g_b^2,x^2(0)\}+\hat{k}^2(0)/\Gamma$ and the bound follows from the inequality $\sqrt{a^2+b^2}\le |a|+|b|$. 
Since $r$ is bounded, this is also true of $q$.}

\textcolor{black}{If $g_{q=r}(t)$ does not vanish identically, $\hat{k}$ must grow without bounds. In this case, $x(t),\dot{x}(t)\rightarrow 0$ (due to the boundedness of $g$ and its partial derivatives with respect to $q$ and $t$) and, consequently, that $q(t)\rightarrow r(t)$ and $u(t)\rightarrow-g_{q=r}(t)$ as $t\rightarrow\infty$, independently of initial conditions. The same conclusions follow if $g_{q=r}(t)\equiv 0$ although, in this case, $\hat{k}$ saturates at some finite positive value.}

\begin{figure}[ht]
	\centering
	\includegraphics[width=.48\textwidth]{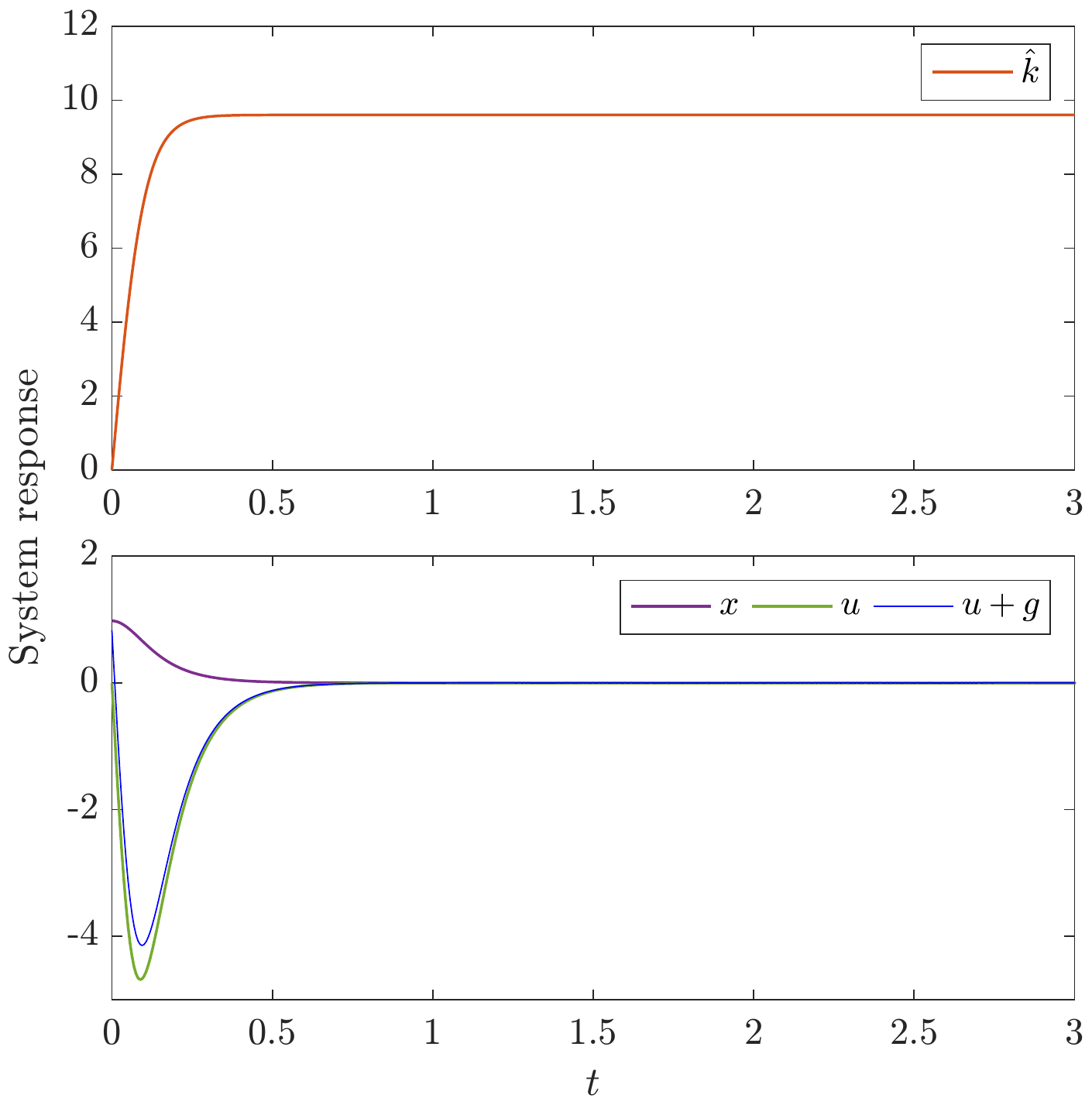}
	\caption{Time histories of $\hat{k}(t)$, $x(t)$, $u(t)$, and $u(t)+g(t)$ of the system defined in \eqref{eq:scalarSys1} at $\omega=1$ under the proposed model-reference adaptive control strategy, where $r(t)\approx q^*(t)$ given by the up-to-fifth-harmonic approximation in \eqref{eq:qstarapprox}.}
	\label{fig:MRAC2Dim1qstar}
\end{figure}

\begin{figure}[ht]
	\centering
	\includegraphics[width=.48\textwidth]{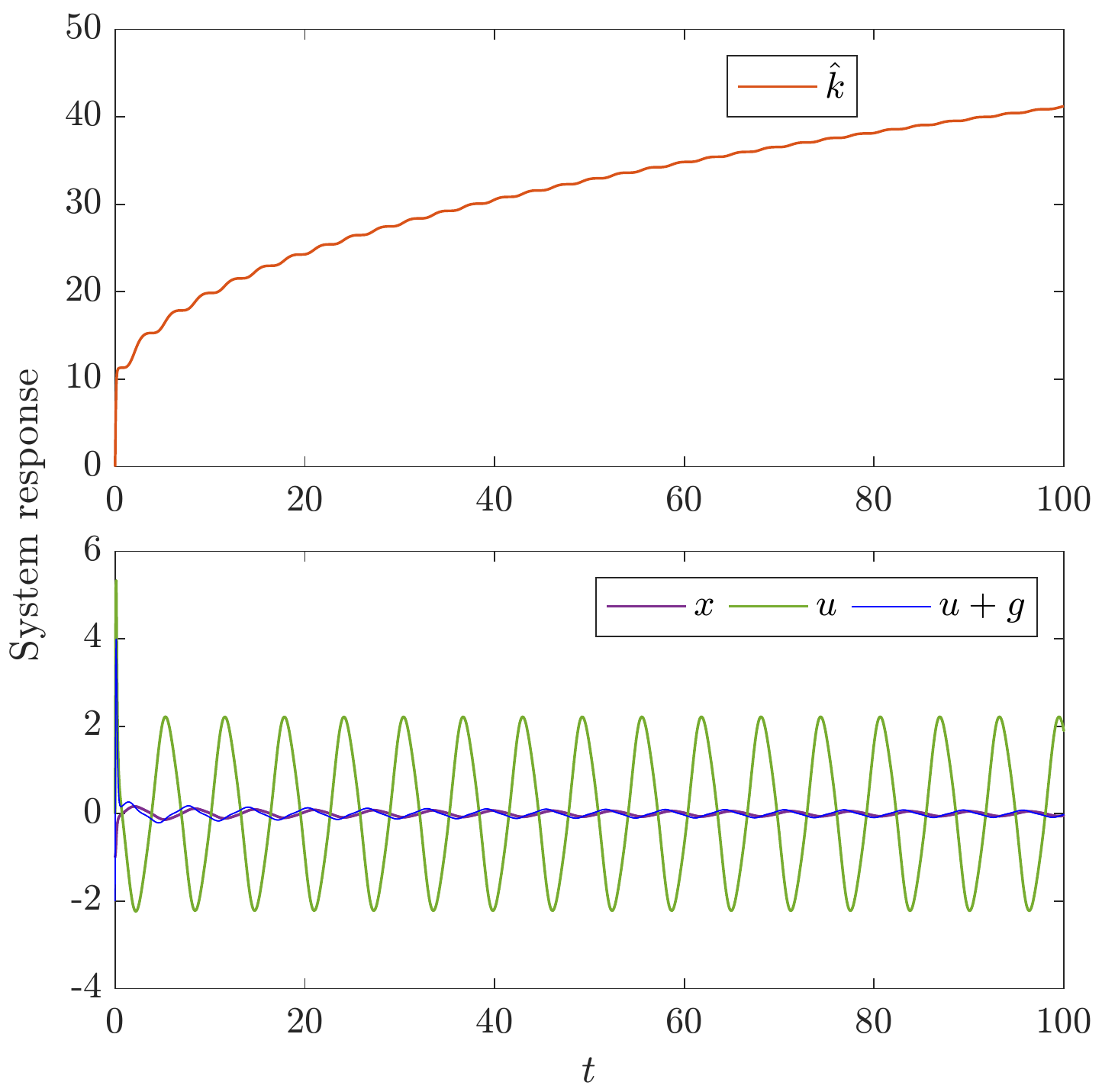}
	\caption{Time histories of $\hat{k}(t)$, $x(t)$, $u(t)$, and $u(t)+g(t)$ of the system defined in \eqref{eq:scalarSys1} at $\omega=1$ under the proposed model-reference adaptive control strategy, where $r(t)=\cos(t)+\sin(t)$.}
	\label{fig:MRAC2Dim1Nqstar}
\end{figure}

These predictions are confirmed by the results of numerical simulations shown in Figs.~\ref{fig:MRAC2Dim1qstar} and \ref{fig:MRAC2Dim1Nqstar} obtained with $r(t)$ given by the right-hand side of \eqref{eq:qstarapprox} and $r(t)=\cos(t)+\sin(t)$, respectively. In each case, $\Gamma=100$ and $q(0)=\hat{k}(0)=0$.

The observations for the example system \eqref{eq:scalarSys1} generalize to a  \textcolor{black}{scalar system of the form 
\begin{align}
\label{eq:NL2SysNctrl}
\dot{q} = aq + b\left(kq+ f\left(t, q\right) + \sigma+u\right),
\end{align}
where $a<0$ and $b\ne 0$} are known constants, $k$ is an unknown constant, and $f(t,q)$ represents a uniformly bounded, unmodeled nonlinearity with uniformly bounded first-order partial derivatives. \textcolor{black}{This reduces to the example system in \eqref{eq:scalarSys1} when $a=-1$, $b=1$, $k=0$, and $f(t,q)=\sin q$.}

With the introduction of a control input $u=-\hat{k}x$ in terms of the deviation $x=q-r$, it follows that\textcolor{black}{
\begin{equation}
    \dot{x}=(a-b\tilde{k})x+g,
\end{equation}
where $\tilde{k}=\hat{k}-k$ and 
\begin{equation}
 g=(a+bk)r + b\left(f\left(t, q\right) + \sigma\right)-\dot{r}.
\end{equation}
If we let
\begin{align}
\label{eq:MRAC2AdaptLawS}
\dot{\hat{k}}&=\Gamma bx^2
\end{align}
and define the Lyapunov function
\begin{equation}
    V=x^2+\frac{1}{\Gamma}\tilde{k}^2,
\end{equation}
it follows by the same argument as for the example that
\begin{align}
\label{eq:MRAC2xbd2}
|x(t)|\leq \max\{-g_b/a,~ |x(0)|\}+\left|\hat{k}(0)-k\right|/\sqrt{\Gamma},
\end{align}
}where $g_b$ denotes the upper bound for $|g|_\infty$. \textcolor{black}{By the uniform boundedness of $f$ and its partial derivatives, we again conclude that $x(t),\dot{x}(t)\rightarrow 0$ and, consequently, $q(t)\rightarrow r(t)$ and $u(t)\rightarrow-g_{q=r}(t)$ as $t\rightarrow\infty$.}

\begin{figure}[h]
  \includegraphics[width=0.48\textwidth]{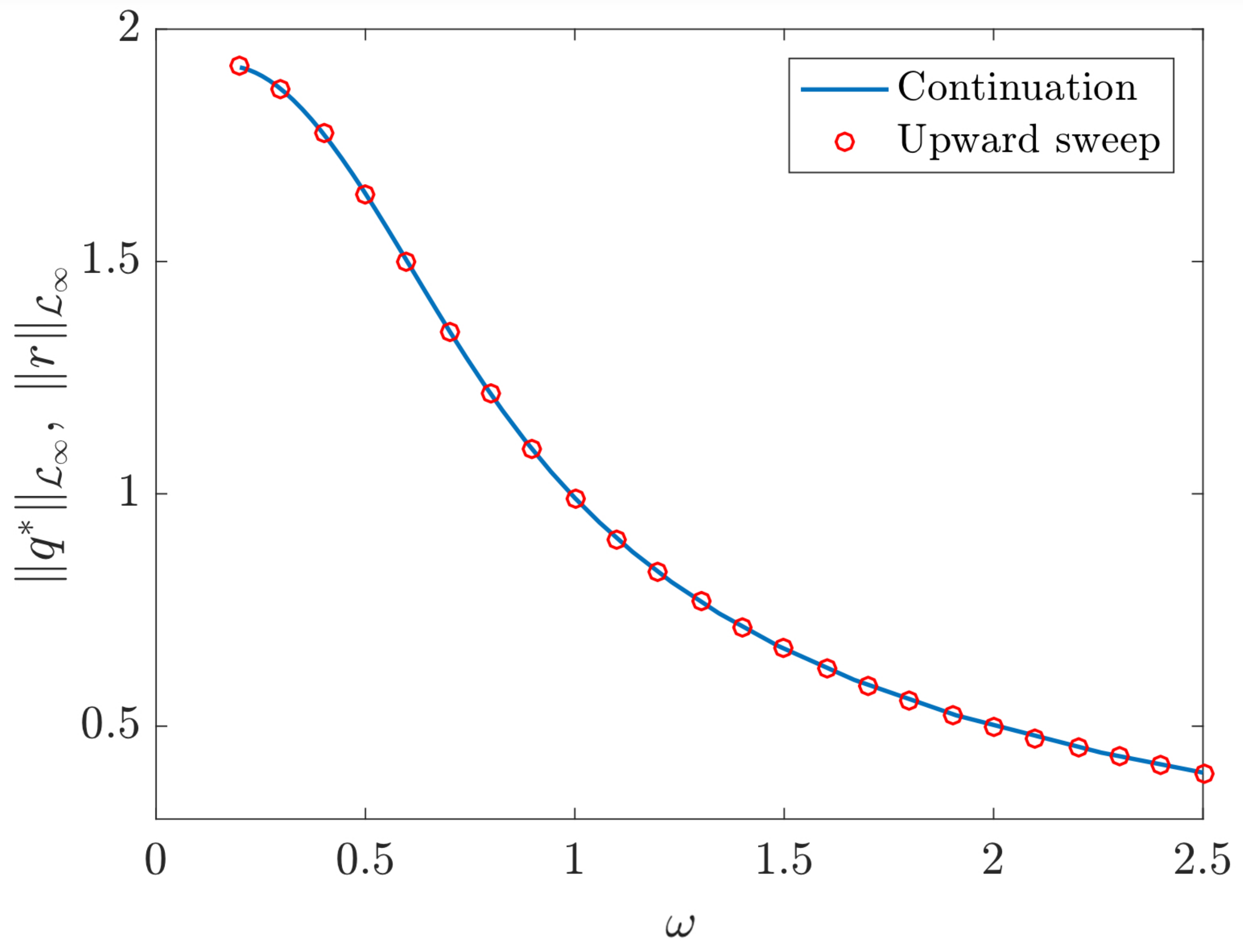}
\caption{A comparison between the branch of periodic solutions of the scalar system defined in \eqref{eq:scalarSys2} obtained using an upward sweep in $\omega$ (represented by $\|q^\ast(\omega)\|_{\mathcal{L}_\infty}$), and that obtained using control-based continuation under the proposed model-reference adaptive control strategy (represented by $\|r(\omega)\|_{\mathcal{L}_\infty}$).}
\label{fig:MRAC2Dim1CBC}
\end{figure}

We proceed to consider the application of the control-based continuation paradigm to the system \eqref{eq:scalarSys1} with $\Gamma=1$. Here, the initial value of the adaptation parameter $\hat{k}(0)$ is set to $0$ in each simulation. Throughout continuation, we approximate $u$ by its truncated Fourier series up to the fifth harmonic and iteratively update the corresponding coefficients of $r$. In each iteration, and when approximating derivatives with respect to the Fourier coefficients of $r$, the Matlab integrator \texttt{ode45} (with relative tolerance $10^{-8}$ and absolute tolerance $10^{-10}$) is used to simulate the closed-loop transient dynamics for $10$ periods, followed by a sampling of the control input during one additional period of simulation. We consider $r$ to have converged to $q^\ast$ when the norm of the Fourier coefficients of the steady-state control input is smaller than $10^{-3}$. The step size along the graph $(\omega,q^\ast)$ is adaptively determined by \textsc{coco} using default settings. The result is shown in Fig.~\ref{fig:MRAC2Dim1CBC}.

Notably, the convergence criterion tolerance is here set orders of magnitude larger than the \textsc{coco} default (which is $10^{-6}$). Indeed, we observe that $x(t)$ remains close to $0$ throughout continuation, resulting in slow dynamics of the adaptation parameter $\hat{k}$ and, consequently, slow rates of convergence of $u(t)$ to $-g(t)$, as required by the control-based continuation paradigm. Although we might be able to improve upon this state of affairs by setting $\hat{k}(0)$ to a larger number, this would likely produce large transient dynamics, including in the control input, making this impractical in physical experiments.

\section{Conclusion}
\label{ch5:concl}

As shown in the previous sections, although adaptive control designs may be proposed for control-based continuation of periodic orbits, significant effort may be required to prove their non-invasiveness and ensure bounded performance and robustness to disturbances, if at all possible. Nevertheless, the benefits over non-adaptive control designs may be equally significant, especially in examples where the linearized dynamics near the sought periodic orbits vary greatly over the parameter range of interest.

It may be reasonably argued that a great number of challenges must be overcome in order to fully realize the potential of control-based continuation already without considering the complexity of self-tuning feedback control designs. Indeed, recent work by Renson, Barton, Sieber, and their collaborators has explored the merger of control-based continuation techniques with data-based approaches such as Gaussian process regression~\cite{Renson202183,Renson20192811} for estimating the local manifold geometry and adaptive filters~\cite{Abeloos20213793} that update the reference input on the fly. The development of the \textsc{continex} toolbox~\cite{Schilder2015251} for \textsc{coco} also highlighted the many challenges associated with continuation in the presence of measurement noise. It would be worthwhile to consider if algorithms for data-based adaptive accommodation of noise and uncertainty along families of periodic orbits could be co-designed with the adaptive feedback control used to locate individual orbits in order to improve overall performance.

For those inclined to explore more general classes of nonlinear problems, for example multi-dimensional systems with unmodeled nonlinearities, we refer to preliminary work described in Chapter 5 of the first author's doctoral dissertation~\cite{Li2019}. There, promising numerical results and guaranteed boundedness of the closed-loop system response \textcolor{black}{(for some proposed choices of adaptive feedback designs)} do not compensate for the lack of proofs of the existence of asymptotically stable limit cycle dynamics, nor the independence of the steady-state control input from initial conditions. Much work remains to be done.

\bibliographystyle{spmpsci}      
\bibliography{nonladapCBC.bib} 

\section*{Statements and Declarations}

\paragraph{Funding}
This work is supported by Agriculture and Food Research Initiative Competitive Grant no. 2014-67021-22109 from the USDA National Institute of Food and Agriculture. Part of the editing of this paper was performed while the second author served at the National Science Foundation. Any opinion, findings, and conclusions or recommendations expressed in this material are those of the authors and do not necessarily reflect the views of the National Science Foundation.

\paragraph{Competing Interests}
The authors have no relevant financial or non-financial interests to disclose.

\paragraph{Author Contributions}
The authors contributed equally to the conception and design of this research, and to the writing of the manuscript. Implementation of algorithms in code and generation of numerical results was performed by Yang Li. Both authors read and approved the final manuscript.

\paragraph{Data Availability}
Datasets generated and analyzed during this study are available upon request from the authors. Matlab scripts sufficient to generate this data will be posted to an open-source archive.

\end{document}